\documentclass[11pt, oneside]{article}   	
\usepackage{geometry}                		
\usepackage{graphicx}				
				
\usepackage{amssymb}

\title{Geometry in the twentieth century: A return to Euclid---The work of Herbert Busemann}

\author{Athanase Papadopoulos\footnote{\emph{Institut de Recherche Mathématique Avancée} (IRMA) and
\emph{Centre de Recherche et d'Expérimentation sur l'Acte Artistique}  (ITI CREAA). Address:
Université de Strasbourg and CNRS,
7 rue René Descartes,
67084 Strasbourg Cedex France;
email: athanase.papadopoulos@math.unistra.fr}}
 \date{\today}							

\begin{document}
\maketitle


 \begin{abstract}
 This is a point of view of the work of Herbert Busemann (1905-1994), seen as a return to the geometry of Ancient Greece. The importance of this work, its recognition and its relation with other works are discussed. The final version of this paper will appear 
 in the Handbook of the History and Philosophy of Mathematical Practice, ed. Bharath Sriraman, Springer, 2024.
 
 \end{abstract}


\section{Introduction}

The twentieth century saw the spectacular development of several fields in mathematics: algebraic topology, low-dimensional manifolds, algebraic geometry,  number theory, dynamical systems, mathematical physics, to name a few. Among these fields, there is one which, in my opinion, is not sufficiently emphasized: metric geometry. In writing this, I am thinking of the work of Herbert Busemann,\index{Busemann, Herbert} whose collected papers have  been published recently in Springer's Collected Works Series  \cite{Busemann-Works-I, Busemann-Works-II}. Busemann's work goes against the tide, and I think that  it is still under-evaluated. This is why I would like to highlight it here. 

The path opened up by Busemann is in a sense a return to Euclid's geometry. I will explain in what sense, and in fact this is the theme of the whole present chapter. Let me first give some indications on what I mean by this return.

To start with, I would like to make it clear that in this chapter, when I  talk about Euclid's geometry, I am not necessarily referring to the propositions found in the \emph{Elements}, but rather to a geometry inspired by them and which is close to them in several respects. To explain this, let me start by recalling that an important  (maybe the most important) part of the geometric books of Euclid's \emph{Elements} is an exploration of the properties of triangles. Likewise, a significant part of Busemann's geometry includes  a study of the geometry of triangles on surfaces (which are not necessarily the Euclidean plane): such a surface may be the sphere, or the hyperbolic plane, or a convex surface embedded in 3-space, or any other two-dimensional manifold. In particular, classical ``non-Euclidean" geometry, as a study of  properties of non-Euclidean triangles, is part of Euclid's geometry as I intend it here. For example, Menelaus' \emph{Spherics}, which I discuss in another chapter of this Handbook \cite{P-Menelaus}, being a treatise on the intrinsic properties of spherical triangles, is an example of what I call Euclid's geometry.  Higher-dimensional hyperbolic and spherical geometry, in particular the study of non-Euclidean polyhedra, which generalize triangles, are also part of this geometry.

The geometry developed by Busemann is defined axiomatically--with a few number of axioms: these include those of a metric space, to which are added some axioms on the existence of shortest lines joining pairs of points. Thus, the notion of distance between two points is one of the basic notions in that geometry. I recall incidentally that Euclid's \emph{Elements} do not use the notion of distance between two points. From the philosophical point of view, this is important and needs a thorough discussion, but for what concerns us here, this is not an issue. Furthermore, in Busemann's geometry,  considerations of differentiability are excluded. Thus, the methods of differential calculus are not part of it. Consequently, Riemannian geometry is not part of this geometry either.

Why call such a geometry ``Euclidean"? Because it deals with problems similar to those raised by Euclid's \emph{Elements} or motivated by them. I shall start by explaining this in the next section, using a significant example, namely, a work of Busemann arising from a question raised by Leibniz. At the same time, the result I am about to describe will show how non-Euclidean geometry naturally fits into this framework.

\section{Leibniz, Busemann and a question of Euclidean geometry}

Leibniz\index{Leibniz!definition of a plane} addressed the question of defining a plane. His work on this subject is reproduced in his \emph{Mathematical writings} \cite[p. 166]{Leibniz}. The problems raised by such a definition are presented and commented on by Busemann in his article \emph{On Leibniz's definition of planes} \cite{Busemann-Leibniz} and are taken up again in his book\index{Busemann, Herbert!\emph{The geometry of geodesics}} \emph{The geometry of geodesics} \cite[p. 46]{GG}. The question of the definition of a plane can be traced back to Euclid, although, as Busemann remarks, the definitions of surface and line that are given in the \emph{Elements} are difficult to work with. We recall that Euclid,\index{Euclid!definition of a line} in Book I of the \emph{Elements}, defines a \emph{line} as something ``that has length and no width", but the terms length and width are not defined there. Likewise, he defines a \emph{plane}\index{Euclid!definition of a plane} as something ``that has only length and width". Busemann notes that ancient authors already found these definitions objectionable, and that  it is not possible, within the framework of Euclidean geometry, to characterize, based only on these definitions, straight lines among curves and planes among surfaces. Thus, one cannot give an answer to Leibniz's question based only on Euclid's \emph{Elements}.  Archimedes went further than Euclid. In \emph{The sphere and the cylinder}, following a definition of straight line due to Heron, he introduced as an axiom the fact that the straight line\index{Archimedes!axiom for a straight line} is a shortest curve among all curves joining two points. He also developed a method for calculating\index{Archimedes!length of curves}  lengths of curves. All this is recalled in the same passage by Leibniz, who considers that this cannot be used for a definition of a straight line, and proposes to define a straight line as the set of points fixed by the rotation of some body---which assumes that the notion of body precedes that of straight line, and  which gives rise to other difficulties. I refer here to Couturat's comments on Leibniz's work, see Chapter IX, of his treatise \emph{The logic of Leibniz} \cite{Couturat},  titled \emph{Geometric calculation}. Regarding Leibniz's question, Busemann writes in the introduction to his article \cite{Busemann-Leibniz}: 
 \begin{quote} \small
 The well-known unsatisfactory Euclidean definitions of straight lines and planes have evoked early attempts to replace them by better ones. Leibniz proposed \emph{to define a plane as the locus of points which have equal distances from two given points}, and the straight line as locus of points having equal distance from three given points. Since the idea of metric space or anything which could be substituted
for it had not yet been developed, this attempt of Leibniz did not grow into a well-founded theory. The present note tries to give a critical investigation of the implications of Leibniz's definition from the point of view of metric spaces.
\end{quote}

To address such a problem, and similar ones, Busemann introduced the axioms of a G-space\index{Busemann!G-space} (the letter G stands for ``geodesic"). In fact, these axioms are already present in Busemann's doctoral thesis (defended in 1931), and they are the same ones he used throughout the rest of his works. This is the moment to introduce these axioms:

\begin{enumerate}

\item The space is equipped with a distance function $d(x,y)$.

\item  Every bounded infinite set has an accumulation point. (Busemann calls a space satisfying this property \emph{finitely compact}.)

\item For any distincts pair of points $x,z$, there is a point $y$ such that $d(x,y)+d(y,z)=d(x,z)$. (One says that \emph{$y$ lies between $x$ and $z$}.)

\item For every point $p$, there exists a positive real number $\rho_p$ such that for any distinct pair of points $x,y$ in the ball of center $p$ and radius $\rho_p$, there exists a point $z$ such that $y$ is between $x$ and $z$. (This is a property of local extendability.)

\item (Uniqueness of extension of geodesics) For any distinct pair of points $x$ and $y$, if $y$ is between $x$ and $z_1$ and between $x$ and $z_2$, and if $d(y,z_1)=d(y,z_2)$ then $z_1=z_2$.

\end{enumerate}

We note that it is not assumed here that the metric space is a manifold. In fact, it is an important open problem, which is presently solved only under some special extra hypotheses, to know whether a G-space is a manifold; see the reviews by  Andreev  \cite{Andreev} and Berestovskii \cite{Berest} who made important contributions to this problem.

The notion of G-space is used extensively in many of Busemann's writings, e.g. in  his 1942 book  \emph{Metric Methods in Finsler Spaces and in the Foundations of Geometry} \cite{Busemann1942},  in his 1943 paper \emph{On spaces in which points determine a geodesic} \cite{Busemann1943}, and in his 1955 book\index{Busemann, Herbert!\emph{The geometry of geodesics}} \emph{The geometry of geodesics} \cite{GG}.

Busemann and B. B. Phadke (one of his former students) write in \cite[p. 190]{Novel}: ``One of the most attractive features of G-space theory is that it maintains a strong connection with the foundations of geometry."

 With the above axioms in the background, Busemann can address Leibniz's question  of defining a plane. He calls a subset of a G-space ``flat" if,  equipped with the induced metric, it is a G-space. In his paper  \cite{Busemann-Leibniz} and in his book \cite{GG}, \S\S 46 and 47, he considers Leibniz's question in the setting of G-spaces to which a single hypothesis is added; it concerns \emph{bisectors}, that is, loci of points
equidistant from two distinct points. The hypothesis says the following:

\medskip

\emph{Any bisector contains with any two points at least one segment connecting them.}

\medskip

Busemann calls this property \emph{flatness} of the bisectors. It is obviously satisfied in the case of the classical geometries of constant curvature. He proves the converse:

\medskip

\emph{The three classical spaces of constant curvature (Euclidean, spherical and hyperbolic) are the only G-spaces in which all bisectors are flat.}

\medskip

Thus, motivated by Leibniz's question, Busemann deduces a remarkable (probably the simplest) characterization of the three classical geometries.
He writes in \cite[p. 307]{GG}, considering such a characterization of the classical spaces: ``Besides being of considerable historical interest, these results form perhaps the most direct and conceptually simplest access to the elementary geometries, and they have many important applications."
 
 In his book \emph{The geometry of geodesics} \cite[p. 310]{GG}, Busemann shows that with a slight modification of the notion of flatness,\index{Busemann, Herbert!\emph{The geometry of geodesics}} one gets a \emph{local} characterization of the three classical geometries. In fact, he gives a characterization of spaces that are \emph{locally} the above spaces of constant curvature, based on a ``local bisector theorem",  a criterion involving a property of small enough open balls; see Theorem 46.1 of \cite{GG}.  It is interesting that Busemann obtains the fact that the metric he ends up with is either the Euclidean, hyperbolic or spherical by establishing the trigonometric formulae for right triangles in such a space, see  \cite[p. 320]{GG}.
 
 Busemann calls the \emph{Helmholtz--Lie problem}  the\index{Helmholtz--Lie problem} question of characterizing the three elementary geometries. 
 From the above bisector theorem characterizing the three classical geometries, he deduces several corollaries, in particular,  statements concerning motions and the Helmholtz--Lie problem, of which the following very general statement is given \cite[Theorem 48.8]{GG}:

\medskip

 \emph{If in a G-space, for any four points $a,a',b,c$ satisfying $d(a,b)=d(a',b)$ and $d(a,c)=d(a',c)$ there exists a motion leaving $b$ and $c$ fixed and carrying $a$ to $a'$, then the space is either locally Euclidean, hyperbolic or spherical.}

\medskip

Another application of the bisector theorem is the fact, obtained in \cite[Theorem 48.16]{GG}, see also \cite{Busemann1947}, that in dimension two, the three elementary geometries are the only geometries in which the area of a triangle is expressible in terms of its sides. More precisely, Busemann proves the following: 

\medskip

\emph{If every point of a G-surface has a neighborhood in which a notion of area for triangles is defined and has the same value for all
isometric triples, then the surface is locally Euclidean, hyperbolic or spherical.}

\medskip

\section{Busemann's  doctoral dissertation}
In one of his last papers, written with Phadke and titled  \emph{Novel results in the geometry of geodesics} (1993)  \cite{Novel}, 
 Busemann recalls that he already introduced the axioms of a G-space (without the name) in his doctoral dissertation, defended  in 1931.
The history of this dissertation is interesting. The thesis supervisor was Richard Courant. Busemann writes in his autobiographical memoirs that he wrote his thesis ``against" Courant. He recounts the following:\footnote{Unless otherwise stated, Busemann’s quotations are taken from the recording of an  interview with Constance Reid, presumably made on April 22, 1973 and
kept at the library of the University of G\"ottingen. I would like to thank again Herbert Goenner for providing me with a copy of this interview. I presented this recording  in more detail in the biography I wrote on Busemann, included in the first volume of his \emph{Collected works}.  \cite[p. 3-12]{Busemann-Works-I}, see also the article \cite{Notices}.}

\begin{quote}\small
Officially, I took my degree with Courant. This was only
``officially", in the sense that I was really inspired by [Pavel Sergueïevitch] Aleksandrov, who visited G\"ottingen regularly. He gave me
the idea of the subject of the thesis. I wrote it, but of course
he could not be the official reviewer of my thesis, so he was
my co-referee.
I must say that my thesis was partly in protest against
Courant. I went to Courant originally, he gave me something
and it turned out to be much easier than what he thought. I
did it and it became a small paper,\footnote{ Busemann talks here about his paper \cite{B-1930}.} but not enough
for a thesis. Then he gave me something else, I did this too,
but then Kolmogorov came to Göttingen, Courant showed
it to him, and Kolmogorov said: ``Oh, that’s all very fine, but
it’s well known." So I became really mad, and I went away
to Rome. I was angry with Courant, and wrote my thesis
there, on my own.
 
\end{quote}

According to the paper by Busemann and  Phadke  \cite[p. 181]{Novel}, 
Busemann became convinced of the importance of non-Riemannian metrics after he read the beginning of Minkowski's \emph{Geometrie der Zahlen} in 1926. At the same time, he  heard a course on point set topology and learned  about  Fréchet’s concept
of metric spaces.\footnote{The concept of metric space as we intend it today,  was first introduced (without the name) by M. Fréchet in his thesis, defended in 1906 \cite{Frechet}.} We read in that paper:
\begin{quote}\small
 The older generation ridiculed the idea
of using these spaces as a way to obtain results of higher differential geometry. But it turned out that a few simple axioms on distance suffice to obtain many non-trivial results of Riemannian geometry and, in addition, many which are quite inaccessible to the classical methods.
\end{quote}

  Busemann wrote his thesis by himself, while he was visiting Rome, and he sent his work to Courant. The latter responded in a letter dated August 8, 1930  that this
``geometric work" is sufficient for a ``reasonable" dissertation.\footnote{An English translation of Courant's report on Busemann's dissertation is contained in Volume II of Busemann's Selected Works.} Let me say a few words about the results of the thesis.

 The subject of Busemann's thesis is the axiomatic
characterization of Minkowski spaces,\index{Minkowski space} that is, metric spaces that underlie finite-dimensional vector spaces equipped with a norm. The paper 
\emph{Über die Geometrien, in denen die ``Kreise mit unendlichem Radius'' die kürzesten Linien sind} (On the geometries where circles of infinite radius are the shortest lines) is extracted from that thesis. An English translation of this paper is contained in Vol. I of Busemann's \emph{Selected Works} 
\cite{On-the-Geometries}. Here, the term ``circle of infinite radius" means a horocycle, that is, the limit of a family of circles passing through a given point and whose center tends to infinity along a geodesic. 
The notion of horocycle (and its higher-dimensional analogue, horosphere) was first used in hyperbolic geometry and it plays a crucial role there.
In Euclidean geometry, the limit of such a family of circles is a straight line. The main result that Busemann proves in the paper \cite{On-the-Geometries} is that 
a 2-dimensional space where horocycles are straight lines is isometric to a Minkowski space. 

The metric characterization of Minkowski spaces remained one of Busemann's main objects of study. I refer the interested reader to the survey titled \emph{Busemann’s work on the characterization of Minkowski geometries} \cite{Papadopoulos-Minkowski} which appeared in the section containing preliminary essays of Volume II of Busemann's \emph{Selected Works} \cite{Busemann-Works-I}.

 In the next section, I will mention other problems considered by Busemann and results he obtained in metric geometry.

\section{Classical notions: questions and results}
The questions and results I mention in this section may seem unrelated but they all pertain to  classical geometry, and many of them date back to the time of Euclid. They are all topics on which Busemann worked. His methods of proof, making use only of elementary geometry, are often ingenious.

    \begin{enumerate}

 \item  A distance characterization of the three classical geometries. Since we already talked about the characterization of the three constant curvature geometries, 
 I start by giving another local distance characterization of these three geometries \cite[p. 47]{GG}:
 
 \medskip
  \emph{If every point in a G-space has a neighborhood in which every quadruple of point is isometric to a quadruple of points in a Euclidean, hyperbolic or spherical space (the choice of such a space may a priori depend on the point), then the space is locally Euclidean, hyperbolic or spherical.}

\item A theory of spaces with prescribed geodesics, and in particular, a thorough study of Hilbert's Problem IV. Roughly speaking, this problem asks for a determination of all metrics on subsets of Euclidean space whose geodesics are the (Euclidean) straight lines. There are several papers by Busemann on this problem, e.g. \cite{B1939, B-Hilbert-IV}; see also the survey \cite{Papadopoulos-Hilbert}. 
 Busemann was the main promoter of Hilbert’s Problem IV,\index{Busemann, Herbert!Hilbert Problem IV} and his work acted as a catalyst to the solution of that problem that was given by A. V. Pogorelov in 1973. The latter writes in the introduction to his monograph \emph{Hilbert’s fourth problem} \cite{Pogorelov-H}: ``The occasion for the present investigation is a remarkable idea due to Herbert Busemann, which I learned from his report to the International Congress of Mathematicians at Moscow in 1966."

\item A theory of perpendicularity\index{Busemann, Herbert!perpendicularity} in G-spaces and conditions under which this relation is symmetric \cite{CIME}.

\item The question of convexity\index{Busemann, Herbert!convexity of spheres} of metric spheres  (that is, sets of points whose distance to a give point is constant) in G-spaces \cite{BP1980, BP1983, Novel}.

\item A theory of parallelism.\index{Busemann, Herbert!parallelism} In the paper \cite[p. 190]{Novel}, Busemann and Phadke write: ``It is hard to understand that the theory of parallels was almost completely neglected in the classical theory after Lobachevsky's discoveries. We now explain the concepts upon which a beautiful theory of parallels can be built in a very general setting." 
 In his 1933 paper \emph{On spaces with convex spheres and the parallel postulate}  \cite{Bu1933},  under the assumption of convexity of metric spheres and the fact that the parallel postulate is satisfied, Busemann studied conditions under which a space can be mapped
 onto a Euclidean space such that the geodesics are mapped to
straight lines. Here, the fact that the \emph{parallel postulate} holds means that if $g$ is an arbitrary straight line (that is, a geodesic isometric to the real line) and $P$ a point not on $g$, then for every sequence of points $x_n$ on $g$ converging to infinity, the line $Px_n$ converges to a unique geodesic.
Parallelism is also the subject of Chapter III of \emph{The geometry of geodesics} \cite{GG}, titled \emph{Perpendiculars and parallels}, see also  the papers \cite{B1979,  Novel} and  the survey \cite{CIME}.  
A theory of parallelism in spaces with non-positive curvature is developed in  \cite[p. 248]{GG}.

\item Angle bisectors in triangles. In his paper  \emph{Planes with analogues to Euclidean angular bisectors} \cite{B-1974}, Busemann  gives a characterization of Minkowski planes based on a property of  angular bisectors of triangles, namely, the fact that an angle in a triangle divides the opposite side in the ratio of the adjacent sides.

\item A theory of length-preserving maps. Here, Busemann's original aim is to determine all locally finite
length-preserving self-maps of the Euclidean, hyperbolic, and spherical spaces \cite{Length-1964}. The first example that he gives of a length-preserving map is folding a piece of Euclidean paper. Among the results he proves is that any finitely compact space with an intrinsic metric and a pairwise transitive group of motions that possesses a proper length-preserving map has constant curvature.

\item The question of existence of  axes of motions, the theory of axial motions, translation distances, asymptotes to axes and parallelism of axes, see \cite[\S 32]{GG} and \cite{Novel} .

\item A theory of horospheres (which he calls  limit spheres)  \cite{Novel}.

 \item A theory of asymptotes,\index{Busemann, Herbert!asymptotes} and the relation between asymptotes and limit spheres \cite{Bu1933, GG, CIME}. The asymptoticity relation is in general neither symmetric nor transitive. Symmetry and transitivity hold in spaces which Busemann calls \emph{straight Desarguesian} \cite[p. 139]{GG}.

\item A study of the properties of fundamental groups of compact spaces of negative curvature  \cite{CIME}.
 
\item The introduction of a theory of timelike spaces: a metric setting for the theory of general relativity, see \cite{Timelike}

\end{enumerate}

The next notions need special attention, and I include them in separate sections.

 \section{The sign of the curvature, angle excess and area of triangles: Replacing equalities by inequalities}
 
Busemann introduced a notion of  sign of the curvature\index{Busemann, Herbert!sign of the curvature} (negative, non-positive, vanishing, nonnegative and positive) in metric spaces. This theory is developed in his book \emph{The geometry of geodesics}  \cite[p. 236 ff.]{GG} and in other works. It is a generalization of the notion of sign of the curvature in Riemannian manifolds in the sense that it coincides with this notion when the underlying space is a Riemannian manifold.  Let us recall the definitions.   In fact, the theory works  for a class of metric spaces that are more general than G-spaces: for most purposes it suffices that the space is \emph{midpoint convex}\index{space!midpoint convex} in\index{midpoint convex space} the sense that for any two points $x$ and $y$, there exists a point $z$ whose distances from $x$ and from $y$ are both equal to half of the distance from $x$ to $y$.  Such a point $z$ (which is generally not unique) is called a \emph{midpoint} of $x$ and $y$.  A G-space is a special case of a midpoint convex space. A midpoint convex space which is complete is geodesic (any two points can be joined by a geodesic).

  Busemann says that a midpoint convex metric space $(M,d)$ is  \emph{non-positively curved}\index{non-positive curvature!Busemann}\index{Busemann!non-positive curvature} if for any three points $a,b,c$ in $M$, if $a'$ is a midpoint of $a$ and $c$ and  $b'$   a midpoint of $b$ and $c$, then
\[d(a',b')\leq 1/2(d(a,b)).
\]

There are variations of this definition, namely, the space is called \emph{negatively curved, non-negatively curved,\index{Busemann!positive curvature}  or\index{Busemann!zero curvature} positively curved}\index{Busemann!non-negative curvature} if \index{Busemann!negative curvature} the\index{Busemann!zero curvature} above inequality is replaced by a strict inequality or by the reverse large or strict inequality respectively. If the inequality is replaced by an equality, the space is said to be \emph{zero curved}. These notions are used today under the name \emph{Busemann non-positive (resp. non-negative, positive, negative, zero) curvature}.

The inequality expressing  negative, non-positive, non-negative, positive curvature 
replaces an equality valid in Euclidean geometry by an inequality.
Obtaining theorems in curved spaces from analogous theorems in Euclidean geometry by replacing equalities by inequalities, is known since antiquity.  This is highlighted in the article \cite{P-Menelaus} contained in the present Handbook  on the work of Menelaus of Alexandria (2nd c. AD). One proposition of 
Menelaus' \emph{Spherics}  says that the sphere is negatively curved in the above sense (see Proposition 27  of  \cite{RP2}), but there many other propositions involving this notion. 
In the \emph{Geometry of geodesics}, Busemann proves that a G-space has non-positive curvature in the above sense if and only if every point has a neighborhood in which all the triangles of sides $\alpha,\beta,\gamma$ and angles $A,B,C$ opposite to these sides satisfy the \emph{cosine inequality} 
\[ \gamma^2\geq \alpha^2+\beta^2-2\alpha\beta \cos C
,\]
that is, if the usual Euclidean cosine formula holds, but where the equality has been replaced by an inequality 
(this is Proposition 41.3 of \cite{GG}).
In the same book, Busemann studies several other metric notions of curvature, see in particular Chapter 5. One of them uses the notion of angle excess of triangles in G-spaces\index{G-space}, see \cite[p. 279ff]{GG}. He also gives several statements that make relations between the theory of parallels and that of angle excess. Let me recall here that the fact that the angle excess of a spherical triangle is positive is also a result of Menelaus, see \cite[Proposition 12]{RP2}. In hyperbolic geometry, angle excess is negative; this was known to J. H. Lambert,\index{Lambert, Johann Heinrich}  for whom the existence of hyperbolic geometry was hypothetical, see \S 81 of his \emph{Theory of parallels} \cite{PT}, in which he defines the area of a triangle in the hyperbolic plane as the \emph{angle defect} of that triangle. Building a theory of angle excess in G-spaces is one of the themes of Chapter V of Busemann's \emph{Geometry of geodesics} \cite{GG}, titled \emph{The influence of the sign of the curvature on geodesics}. 

Talking about area, let me quote another result of Busemann giving a characterization of Minkowski spaces in terms of a property of areas of its triangles. This is Theorem 51.5 of \cite{GG}. The theorem concerns  G-spaces that are \emph{straight}, that is, in which the property of local extendability in the axioms of a G-space is replaced by a global one.
The result says that in a straight space, the fact that the area of any 
 triangle $ABC$ depends only on the distances $BC$ 
and the distance from $A$ to the segment $BC$ is equivalent to the fact that the space is a
Minkowski plane with symmetric perpendicularity.

\section{Convexity theory}
 
 Among the introductory essays\index{Busemann, Herbert!convexity} included in Volume II of the \emph{Selected Works}, there is one titled \emph{Busemann and Feller on curvature properties of convex surfaces} \cite{ACP}, written by Annette A'Campo and the author of the present chapter.   It  is   a commentary on the paper \emph{Krümmungseigenschaften konvexer Flächen} (Curvature properties of convex surfaces), written in 1935, by Busemann and Feller \cite{Bus-F}.

 This work of Busemann and Feller is strongly influenced by Hilbert's\index{Hilbert, David (1862--1943)} vision of mathematics. Both Busemann\index{Busemann, Herbert} and Feller\index{Feller, William (1906--1970)} obtained their PhD in G\"ottingen, which, in those times, was one of the two most important mathematical centers of the world (the other one being Paris). Although their thesis advisor was Courant, both Busemann and Feller were shaped by Hilbert's\index{Hilbert, David (1862--1943)} view on the foundations of geometry. Hilbert retired in 1930 but continued to give a course on the philosophy of mathematics. Feller\index{Feller, William (1906--1970)} obtained his PhD in 1927, and Busemann in 1931. Before talking about their work, let me say a few words on some work of Hilbert related to the same topic.

 In 1932, together with Stephan Cohn-Vossen,\footnote{For a comprehensive review on Cohn-Vossen's work, we refer the reader to the article \cite{A-Cohn} by A. D. Aleksandrov.}\index{Cohn-Vossen, Stephan (1902--1936)} Hilbert\index{Hilbert, David (1862--1943)} published his remarkable book \emph{Geometry and the Imagination}\footnote{The original German title is: Anschauliche Geometrie.} \cite{HilbertCohnVossen}  which is based on Hilbert's lectures held in 1920/21. The book contains a chapter on the geometry of surfaces, but the point of view is purely geometric, unlike the one of Euler, Monge and most of the other mathematicians who followed, who studied surfaces and curves using analytic geometry. The methods used in the chapter  on differential geometry of Hilbert and Cohn-Vossen's book  constitute a return to the methods of the Greeks. For instance, their description of the osculating circle at a point of a plane curve starts with: 
`` We draw a circle through $P_1$ and two neighboring points on the curve. If we let the two neighboring points approach $P_1$, the corresponding circle converges to a limit position [\ldots] This limit circle is called the osculating circle\index{osculating circle} at $P_1$." No differentiation is involved.
The discussion then is generalized to space curves.  The osculating plane of a space curve at a point  is constructed 
as the limit of the plane through the tangent at $P$ and a neighboring point that approaches $P$   \cite[p.\,~158]{HilbertCohnVossen}. After a study of intrinsic properties that characterize the sphere, Hilbert\index{Hilbert, David (1862--1943)} and Cohn-Vossen\index{Cohn-Vossen, Stephan (1902--1936)}
In \cite{HilbertCohnVossen}, \S 32 study convex surfaces sharing the same properties.

Between the work of Hilbert and Cohn-Vossen  and that of Busemann and Feller on convexity, I mention the work of Bonnesen and Fenchel on the same topic. 
They wrote a book titled \emph{Theorie der konvexen K\"orper} (The theory of convex bodies) which was published in 1934, the year Werner Fenchel, who also belonged to the Göttingen school, moved to Denmark. This book, which is  one of the founding texts of modern convexity theory,\index{convexity} is also written in the tradition of Hilbert's foundations of geometry.  Section 17 starts with the following: ``Far-reaching assertions can be made about the curvature relations of convex curves without making any assumption of differentiability" \cite[p.\, 153]{BF}.  Indeed, the authors, in that section and the next one,  develop a geometric theory of curvature\index{curvature!surfaces} of convex curves and surfaces without differentiability conditions.  
 
  Busemann's works on convexity is in the same style.
 His first published work on the subject is the joint paper  \cite{Bus-F} with Feller.\footnote{William Feller (1906-1970) studied in Götttingen, obtaining his PhD with Courant, like Busemann. Feller obtained his PhD in 1927, and  Busemann in 1931. They both emigrated to the USA after the Nazis formed a coalition in Germany. They wrote seven joint papers, published between 1934
and 1945.  The reader interested in this turbulent period of the history of Europe in relation with mathematics may refer to Feller\index{Feller, William (1906--1970)} and Busemann's\index{Busemann, Herbert} biographies contained in their Selected Works editions, \cite[Vol. 1]{FellerWorks} and  \cite{Busemann-Works-I}.} In this paper, the authors study curvature properties of
convex surfaces  in 3-space.  Notions like tangent line, tangent plane, osculating
circle, radius of curvature, etc. are studied using their 
original geometrical meaning rather than differential calculus. Busemann and Feller prove that a convex surface has almost everywhere
second-order derivatives. (It was only known before that such a surface is almost
everywhere differentiable, cf. \cite{Reidemeister}.)  Among the other results, they obtain that there is a
 well-defined nearest point projection of curves lying in
the complement of a convex body bounded by a convex surface onto this
convex surface which is length non-increasing. They also study geodesics
 on convex surfaces. 
  They describe their goal in the introduction as follows: 
``The aim of this paper is to develop the elements of differential geometry of convex surfaces without making the usual regularity assumptions. Hereby necessarily, purely geometric considerations often replace the usual analytic arguments."\footnote{Cf. the English translation in \cite[p.~235]{Busemann-Works-I}.}

 They show that without any differentiability
hypothesis, one can define a notion of
curvature and obtain the analogues of  the classical theorems of Euler, Meusnier and others on the differential
geometry of surfaces. Their arguments are synthetic.  
The theorems of Meusnier and Euler are recurrent in Busemann’s work.
In a paper he wrote 15 years later \cite{Busemann1950}, Busemann, introduced a new
notion of curvature in Minkowski spaces, and to justify the relevance of this notion, he showed that analogues
of the theorems of Meusnier and Euler hold in this new setting as well.  
The study of the differential geometry of curves and surfaces, based on
the notion of tangent, normal, order of contact, etc. of plane curves finds its
origin in the work of Archimedes and Apollonius (third and second centuries
B.C.) and other Greek geometers.

This work of Busemann and Feller was the beginning of a lifelong program
conducted by Busemann on the study of convexity theory using a minimal amount of differential geometry or differential calculus.
He worked later on various geometric aspects of this subject, including 
 volume, Brunn--Minkowski theory and convexity in general metric
spaces.

The subject of convexity dates back to Greek antiquity. Busemann, in his study of convex curves in G-spaces in his book \emph{The geometry of geodesics} \cite[p. 154]{GG}, after he proves that a convex simple closed curve in 	a 2-dimensional G-space has a supporting line at each point, writes that ``it is remarkable that an exact proof of this proposition in the plane was already given by
Archimedes." On the notion of angle, tangents, and contact in Greek mathematics, I refer the reader to the remarkable editions by R. Rashed, \cite{R1, R2, R3, R4, R5, R6}.

   In 1956, together with C. M. Petty,  Busemann formulated a set of problems on convex bodies, in a paper titled \emph{Problems on convex bodies} \cite{B-Petty}.
Most of these problems are still open, see the commentary by D. Burago \cite{Burago-M} in Volume II of Busemann's \emph{Selected Works}. The latter, in the conclusion to his  article, writes that ``[he believes] that this circle of problems posed by
Busemann opened a can of worms, a new direction or perhaps even a new area of
research which is likely to be very active for many years to come."

In 1958, Busemann published a book titled \emph{Convex surfaces}  \cite{Busemann-convex}. 
In the Introduction, he declares that his purpose is to present a subject, convex surfaces, which ``during the past 25 years has experienced a striking and beautiful development, principally in Russia, but has remained largely unknown, at least in the USA." Indeed, the book is a tribute to the work of A. D. Aleksandrov.\index{Aleksandrov, Alexandr Danilovich} The following paragraph gives us some information on how Busemann presented  his book project to the publisher:

A letter dated July 23, 1956, from the editor-in-chief of Interscience Publishers Inc., addressed to Busemann, starts with:\footnote{The letter is part of Courant's correspondence kept at the Elmer Holmes Bobst Library in New York. Carbon copies of the letter were sent to Courant and to Bers.} ``It is with very special pleasure indeed that we learn from Dr. Stocker of your willingness to present in our Tracts series an account of important recent development of the work of certain Russian geometers." The letter refers to the work of A. D. Aleksandrov and the school he founded in Russia. I am not going to attempt any review of the results of Aleksandrov and his school because this would lead us much too far, but I shall mention some of this work  later in this chapter. It is interesting to know that Busemann's book on convexity which was initially planned to be a review of the work of the Russian school on this subject, was eventually translated into Russian in 1964, see \cite{Busemann-convex}. The book contains several important results on convex bodies.
Among these results, let me mention a proof of a theorem in convexity theory that was independently obtained by Alexandrov\index{Aleksandrov, Alexandr Danilovich}  \cite{Alex1937} and Fenchel\index{Fenchel, Werner}  and Jessen \cite{FJ1938}, stating that if two convex bodies in Euclidean $n$-space have equal $p$-th order area measure (we shall not go into the definition here) for some $p$ satisfying $1\leq p\leq n-1$, then the two bodies are the same up to a translation.

\section{Busemann's work recognition}
The recognition of of Busemann's work came from the Soviet Union. In 1985, he was awarded the
Lobachevsky Prize  ``for his innovative book \emph{The Geometry of geodesics}". The book was written 30
years before. The Lobachevsky Prize (or Lobachevsky medal, depending on the epoch and on the institution which delivered it), is still a prestigious prize. Its first recipients were Sophus
Lie in 1897, Wilhelm Killing in 1900 and David Hilbert in 1903. A.
D. Aleksandrov obtained it in 1951. The list of the more recent recipients  includes
Weyl, Pontryagin, H. Hopf, P. S. Alexandroff, Efimov, Kolmogorov, Hirzebruch,
Arnol’d, Margulis, Gromov and Chern.

Busemann's work was not appreciated for its true value in the West. In the 3-volume set titled \emph{A Century of Mathematics in America} \cite{Duren} published in 1988 by the AMS, there is no mention of Busemann's work. 
In a letter to 
V. Pambuccian, dated Nov. 26, 1989 and published in Volume II of Busemann's \emph{Selected Works} \cite[p. 64]{Busemann-Works-II}, Busemann writes: ``\ldots Your complaint regarding the neglect of synthetic geometry in the USA is
unfortunately justified.
\emph{The Geometry of Geodesics}  (GG) and \emph{Recent synthetic Differential Geometry} are still available,
however for ridiculously high prices. \ldots The absence of interest in the USA is recognizable from the very fact that every year there are approximately
twice as many copies sold abroad than in the USA. And the factor two is probably lower than
the real one, for three of my books were translated in Russian, and appeared in well printed, large, very
cheap editions, e.g. GG costs 1.88 rubles."

Busemann's work started to be recognized in the West only in the 1980s, when metric geometry was revived  by M. Gromov, who introduced  the methods of synthetic
global geometry. In fact, metric geometry became fashionable in France after the latter  started giving courses on this topic in Paris at the end of of the 1970s. His 1979 course was published in several editions, see \cite{Gromov-Structures}.  Igor Belegradek, in his \emph{Math Reviews} review of the first English edition  of this book (1991), writes that it is ``one of the most influential books in geometry in the last twenty years". This coincided with the last years of Busemann's career.   Gromov is an heir of the Russian school of metric geometry founded by A. D. Aleksandrov. I will talk about the latter in the next section.

\section{On the work of A. D. Aleksandrov}

 A. D. Aleksandrov founded an important school
on metric geometry, with a large number of  students and collaborators.
Like Busemann, he was interested in the most basic notions of geometry. Classical problems of convexity, isoperimetry and
isoepiphany, on which he and the school he founded were working, became at the forefront of research. For him, like for Busemann,
 classical synthetic
geometry and the old geometric problems that originate in Greek Antiquity were more important than Riemannian geometry based
on linear algebra and tensor calculus. In some sense,
his work was a return to Euclid and Archimedes.

Like Busemann,\index{Busemann, Herbert} Aleksandrov's\index{Aleksandrov, Alexandr Danilovich} interest abides in the most basic notions of geometry. 
Aleksandrov's  work was perpetuated by his students and his student's students. Yurii Burago, Yurii Borisov, Ilya Backelman, Alexei Pogrelov, Victor Zalgaller, Grigorii Perelman, Yurii Reshetnyak, Anton Petrunin,  and Valerii Berestovskii are among his descendents---there are many others. Aleksandrov was an invited speaker at the 1960 ICM, where  he delivered a talk titled \emph{Modern develompent of surface theory} \cite{Alex-Modern} in which he presented an outline of his work together with that of his students, including the metric condition of the boundedness of curvature of triangles in terms of angle excess leading to  the metric definition of curvature, and a variety of  notions and results obtained under the condition of bounded curvature, including a definition of area, a study of convexity, embedding problems, parallel translation, and several others.
 We saw that in the paper  \cite{Bus-F},  Busemann and Feller proved that
a convex surface in 3-dimensional space has a second order differential at almost every point, and
that the Dupin indicatrix has the usual form. In 1939, A. D. Aleksandrov extended
this result to an arbitrary convex hypersurface in $n$-dimensional Euclidean space.   I refer the reader to Aleksandrov's\index{Aleksandrov, Alexandr Danilovich}  book \emph{Intrinsic geometry of convex surfaces} \cite{Aleksandrov-II}, in which he provided proofs of the theorems of Euler, Meusnier\index{Meusnier Theorem}\index{Theorem!Meusnier} and Rodrigues\index{Rodrigues formula}\index{formula!Rodrigues} in the general setting of convex surfaces, attributing these generalizations to Busemann and Feller. In the same book,  he proved that any convex  two-dimensional manifold is non-negatively curved in the sense he introduced and which we recalled above, and that every non-negatively curved two-dimensional manifold homeomorphic to the sphere is isometric to some convex surface.

We mentioned Busemann's work on ``timelike spaces", which provides a geometric setting to general relativity.  Aleksandrov with his students developed a similar theory, which they called \emph{chronogeometry}, see the review in \cite{Papa-Chronogeometry} and the references there.

The metric notions of curvature introduced by Aleksandrov and by Busemann turned out to be at the bases of very important developments in geometry, and the terms ``Busemann geometry" and ``Aleksandrov geometry" refer to these works. Two AMS Mathematics Subject Classification codes that carry the names of Busemann and Aleksandrov geometries\footnote{These are: 53C70: Direct methods (G-spaces of Busemann, etc.) and 53C45: Global surface theory (convex surfaces à la A. D. Aleksandrov).} and there is a multitude of theorems and research problems that rely on their intrinsic definitions of the sign of the curvature.

In \cite{KNR} Kutateladze, Novikov, and Reshetnyak write:
\begin{quote}\small
 It would be impossible to understand Aleksandrov’s philosophy without turning
to the roots of the science he loved best. In 1981 he wrote: `The spirit of modern
mathematics lies in the tendency to return to the Greeks'. His favorite motto was:
``Back to Euclid!" Geometry was a part of the culture of the ancient world.
\end{quote}

 They also write :
\begin{quote}\small
 Aleksandrov steered mathematics towards the synthetic geometry of the ancients
in a much subtler way than is now recognized. What is meant is not just a transition
from smooth local geometry to geometry in the large without restrictions
on differentiability. By enriching differential geometry with the tools of functional
analysis and measure theory, he strove to regain the grandeur of mathematics in the
times of Euclid. And by synthesizing geometry with other branches of 20th-century
mathematics, he aspired to the ancient ideal of a unified science.
  \end{quote}
 
See also the article by Katateladze  \cite{Kutateladze}, with the significant title 
\emph{Aleksandrov of Ancient Hellas}. Likewise, Zalgaller, in his historical paper \cite{Zalgaller} in which he describes the lively Leningrad geometry  seminar led by Alexandrov, writes about the latter: 
``When talking about methodology, he often pronounced the slogan: `Retreat to Euclid!', fighting for visual methods."

%
%
%
%
%
%
%
%
%

\section{In guise of a conclusion}

Since his doctoral dissertation and until the end of his life, Busemann used to formulate problems and work on them without relying on the trends that were
fashionable in his time. In an article that appeared in the Los Angeles
Times on June 14, 1985, \cite{Dembart}, on the occasion of the attribution of the
Lobachevsky Prize to Busemann, the author writes:
``Few mathematicians ever make it into public consciousness,
but Busemann has had a hard time even within his own
field, in part, at least, because he never worked on trendy
problems and never followed the crowd." The same journalist quotes Busemann saying: ``If I have a merit, it is that
I am not influenced by what other people do."

A geometry freed from the dominion of differential calculus kept growing since the last quarter of the twentieth century, reinvigorated in the works of Thurston\index{Thurston, William P. (1946--2012)} whose methods constitute a return to the synthetic methods of non-Euclidean geometry, and of Gromov,\index{Gromov, Mikhail}  closer in spirit to the Busemann and Aleksandrov approaches. Problems on area and volume,  isoperimetry and convexity, became the catalyzers of mathematical research.

 Finally, I started this article by mentioning Leibniz. I would like to end it by mentioning another mathematician-philosopher, René Thom.

Thom was one of the main founders of modern topology. His work on cobordism,  for which he was awarded the Fields medal in 1958, inaugurated the combination of the methods of differential topology with those of algebraic topology.
 He writes in \cite[p. 10]{Predire} that in his youth, ``[he] had a certain taste for Euclidean geometry which immediately attracted [him] very much."
He kept this taste for Euclidean geometry---the ``old Euclidean geometry", as he used to call it, with its constructions, criteria for triangle equality, etc. all his life, even though, in the mathematical milieu in which he was to evolve, it was against the grain.  Having disappeared from the school and university curricula. Euclidean geometry became a field not only unknown, but, what is worse, he said, poorly understood. He writes in  \cite{Thom-Logos}:  ``I have kept a preference for Euclidean geometry that my colleagues and fellow mathematicians do not forgive me."

  Busemann's work is still poorly known and  
greatly underrated. I find no good reason for that, except that this works lacks publicity, which, it seems, Busemann himself was reluctant to do.  In his letter to V. Pambuccian, dated Nov. 26, 1989, which I quoted above, he writes:
``The reason for the neglect of synthetic geometry appears to me to lie first in the fact that every professor
cares about seeing his ideas live on via his students. To admit that his methods use to a large extent
unnecessary hypotheses is inconvenient and is simply swept under the rug. In fact, I do not know, except
for some of my students, anybody in the USA who has actually read GG."

\bigskip

\noindent{Acknowledgement} This work is supported by the lnterdisciplinary Thematic lnstitute CREAA, as part of the ITI 2021-2028 program of the Université de Strasbourg, the CNRS, and the Inserm (funded by IdEx Unistra ANR-10-IDEX-0002, and by SFRI-STRAT’US ANR-20-SFRI-0012 under the French Investments for the Future Program).


\begin{thebibliography}{00}  
  
  \bibitem{ACP} A. A'Campo Neuen and A. Papadopoulos, Busemann and Feller on curvature properties of convex surfaces, Busemann and Feller on curvature properties of convex surfaces, In: Herbert Busemann, Selected works, II, ed. A. Papadopoulos,
Springer, Cham, 2018
p. 67-88.

   
 \bibitem{Alex1937} A. D. Aleksandrov, Zur Theorie der gemischten Volumina von konvexen Körpern. II. (russ.) Mat. Sbornik N. S.  2, (1937), p. 1205-1238.


\bibitem{Aleksandrov1939} A. D. Aleksandrov, Almost  everywhere  existence  of  the  second  differential  of  a  convex
function  and  some  properties  of  convex  surfaces  connected  with  it, Leningrad  State Univ. Annals [Uchenye Zapiski] Math. Ser. 6
(1939), p. 3-35.




    \bibitem{A-Cohn} A. D. Aleksandrov, On the works of S. E. Cohn-Vossen, Uspehi Matem. Nauk (N.S.) 2, (1947), no. 3 (19), p. 107-14. English translation by  A. Iacob.

 
 \bibitem{Aleksandrov-II} A. D. Aleksandrov, Intrinsic geometry of convex surfaces, ed. S. S. Kutateladze,  transl. by S. Vakhrameyev, Boca Raton, FL: Chapman \& Hall/CRC Press, 2005.
 
  \bibitem{Alex-Modern} A. D.
Aleksandrov, Modern development of surface theory.
Proc. Int. Congr. Math. 1958, 3-18 (1960).

  
  
 \bibitem{Aleksandrov-III} A. D. Aleksandrov and V. A. Zalgaller, Intrinsic geometry of surfaces, Academy Sc. USSR, Moscow, 1962. English translation by J. M. Danskin, AMS Translations of Mathematical Monographs, Amercian Mathematical Society, Providence, 1967.
 
 
  \bibitem{Andreev} P. D. Andreev, Busemann's problems on G-spaces,  In  H. Busemann, Selected works, II, ed. A. Papadopoulos,
Springer, Cham, 2018, p. 85-113.

  
\bibitem{Ber1982} V. Berestovskii,  Busemann homogeneous G-spaces.
Dokl. Akad. Nauk SSSR 247 (1979), no. 3, 525–528 


  \bibitem{Ber1979} V. Berestovskii,  Homogeneous Busemann G-spaces.
Sibirsk. Mat. Zh. 23 (1982), no. 2, 3-15, 215.


\bibitem{Berest} V. N. Berestovskii, Busemann's results, ideas, questions and locally compact homogeneous geodesic spaces, In:  H. Busemann, Selected works, II, ed. A. Papadopoulos,
Springer, Cham, 2018, p. 41-84.


\bibitem{Berest-1} V. N. Berestovskii, Herbert Busemann and convexity,  In: H. Busemann, Selected works, 2, ed. A. Papadopoulos,
Springer, Cham, 2018, p. 89-118.


 \bibitem{BF} T. Bonnesen and W. Fenchel, Theorie der konvexen K\"orper, Springer, Berlin, 1934,  English translation: Theory of convex bodies,  BCS associates, Moscow, Idaho, 1987.

\bibitem{Burago-M} D. Burago, On Several Problems Posed by H. Busemann and Related to Geometry
in Normed Spaces. H. Busemann, Selected works, II, ed. A. Papadopoulos,
Springer, Cham, 2018,  p. 119-121.

\bibitem{B-1930} H. Busemann, 
Die Vollständigkeit der Minimalfolgen von Eigenwertproblemen,
Nachrichten Göttingen 1930, 295-307 (1930). 
 
\bibitem{Busemann1950}  H. Busemann, The foundations of Minkowskian geometry. Comment. Math. Helv. 24,
(1950), 156-187





\bibitem{On-the-Geometries} H. Busemann, Über die Geometrien, in denen die ``Kreise mit unendlichem Radius'' die kürzesten Linien sind, 
Math. Ann. 106 (1932), no. 1, 140-160. English translation by A. A'Campo, On the geometries where circles of infinite radius are the shortest lines, in :  Vol. I of \cite{Busemann-Works-I}, p. 155-172.


\bibitem{Bu1933} H. Busemann, Über Räume mit konvexen Kugeln und Parallelenaxiom.  
Nachrichten Göttingen 1933, 116-140 (1933).  English transl. by A. A'Campo, On Spaces With Convex Spheres and the Parallel Postulate,  In  H. Busemann, Selected works, I, ed. A. Papadopoulos,
Springer, Cham, 2018,  p. 173-191.


 
  
  
  \bibitem{B1939} H. Busemann Two-dimensional metric spaces with prescribed geodesics.  
Ann. Math. (2) 40, 129-140 (1939). 


\bibitem{Busemann-Leibniz} H. Busemann, On Leibniz's definition of planes,  Amer. J. Math. 63 (1941), 101-111.




\bibitem{Busemann1942} H. Busemann,   
 Metric Methods in Finsler Spaces and in the Foundations of Geometry.
Annals of Mathematics Studies, no. 8. Princeton University Press, Princeton, N. J., 1942.



\bibitem{Busemann1943} H. Busemann, 
On spaces in which points determine a geodesic.  \emph{Trans. Am. Math. Soc.} 54  (1943), 171-184. 

\bibitem{Busemann1947} H. Busemann, Two-dimensional geometries with elementary areas. Bull. Amer. Math. Soc. 53 (1947), No. 4, 402-407.


\bibitem{Busemann-1948} H. Busemann, Spaces with non-positive curvature, Acta Math. 80 (1948), 259-311.


\bibitem{Busemann-Brunn}  H. Busemann, A Theorem on Convex Bodies of the Brunn--Minkowski type, Proc. Natl. Acad. Sci. USA 35, 27-31 (1949).


\bibitem{Busemann-convex} H. Busemann, Convex surfaces, Interscience tracts in pure and applied mathematics, Interscience Publishers, New York and London, 1958.  Russian translation, Convex surfaces,
Izdat. Nauka, Moscow, 1964.


\bibitem{B-1974}  H. Busemann, Planes with analogues to Euclidean angular bisectors, Math. Scand. 36 (1975), 5-11.
  

\bibitem{GG} H.  Busemann,    The geometry of geodesics.
Academic Press, Inc., New York, N.Y., 1955, reprinted by Dover in 2005.  Russian translation, Gosudarstv. Izdat. Fiz.-Mat. Lit., Moscow, 1962, 504 p. 



 \bibitem{CIME} H. Busemann, The synthetic approach to Finsler spaces in the large,
C.I.M.E. Summer Sch., 23, Geometria del Calcolo delle Variazioni (1961), p. 1-72, Springer, Heidelberg.



\bibitem{Length-1964} H. Busemann, 
 Length-preserving maps
Pacific J. Math. 14 (1964), p. 457-477.



\bibitem{Timelike} H. Busemann, Timelike spaces.  
Diss. Math. 53, 52 p. (1967). 


  \bibitem{B-Synthetic} H. Busemann, 
Recent synthetic differential geometry. Ergebnisse der Mathematik und ihrer Grenzgebiete, Band 54. Berlin-Heidelberg-New York, Springer-Verlag, 1970.   

\bibitem{B-Hilbert-IV} H. Busemann, Problem IV: Desarguesian spaces.  
Math. Dev. Hilbert Probl., Proc. Symp. Pure Math. 28, De Kalb 1974, 131-141 (1976). 



\bibitem{B1979} H. Busemann, Minkowskian geometry, convexity conditions and
the parallel axiom. J. Geom. 12, no. 1, p. 17-33, (1979).


\bibitem{Busemann-Works-I} H.  Busemann,    Selected works, I, ed. A. Papadopoulos,
Springer, Cham, 2018. 


\bibitem{Busemann-Works-II} H.  Busemann,    Selected works, II, ed. A. Papadopoulos,
Springer, Cham, 2018.


\bibitem{Bus-F}  H. Busemann and W. Feller, Krümmungseigenschaften Konvexer Flächen,
Acta Math. 66 (1936), no. 1, 1-47, English translation by A. A’Campo Neuen: Curvature Properties
of Convex Surfaces. Selected Works, Herbert Busemann, Vol. II, p. 235–269..
 
 \bibitem{B-Petty} H. Busemann and C. M. Petty, Problems on convex bodies, Math. Scand. 4 (1956), 88-94.



\bibitem{BP1980} H. Busemann and B. B. Phadke, Nonconvex spheres in G-spaces.  
J. Indian Math. Soc., New Ser. 44, 39-50 (1980). 


\bibitem{BP1983} H. Busemann and B. B.  Phadke,  Peakless and monotone functions on G-spaces
Tsukuba J. Math. 7 (1983), no. 1, 105-135 : 

\bibitem{Novel} H. Busemann and B. B. Phadke,  Novel results in the geometry of geodesics
Adv. Math. 101 (1993), no. 2, 180-219.


\bibitem{Couturat} L. Couturat, La logique de Leibniz : d'après des documents inédits, Paris, Alcan, 1901.
 
\bibitem{Dembart}  L. Dembart, An unsung geometer keeps his own plane, Los Angeles Times,
June 14, 1985.



 
\bibitem{Duren} P. L. Duren, R. Askey, U. C. Merzbach, H. M. Edwards (ed.),  A Century of Mathematics in America,
American Mathematical Soc., 1988.


\bibitem{FellerWorks} Feller, W. Selected papers,  I.
R. L. Schilling, Z. Vondraček, W. A. Woyczynski (ed.), Springer, 
Cham, 2015. 


 
 \bibitem{FJ1938} W. Fenchel and B. Jessen, Mengenfunktionen und konvexe Körper, Danske Vid. Selskab Mat.-Fys.Medd. 16, 3 (1938), p. 1-31. 
 
  
 \bibitem{Frechet} M. Fréchet, Sur quelques points du calcul fonctionnel, Thèse présentée à la Faculté des Sciences de Paris pour obtenir le grade de Docteur ès Sciences, 1906, published in
Rendiconti del Circolo Matematico di Palermo, vol. 22, p. 1-72 (1906).
 
\bibitem{Gromov-Structures} M. Gromov,
Structures métriques pour les variétés riemanniennes, Edited by J. Lafontaine et P. Pansu.
1st ed. Textes Mathématiques, 1. Paris: Cedic/Fernand Nathan. vii, 152 p., 1981. English translation: 
Metric structures for Riemannian and non-Riemannian spaces. Transl.S. M. Bates. With appendices by M. Katz, P. Pansu, and S. Semmes.  3rd printing, Modern Birkhäuser Classics. Basel: Birkhäuser, 2007. 

\bibitem{HilbertCohnVossen} D. Hilbert and S. Cohn-Vossen, Anschauliche Geometrie, Springer, Berlin and Heidelberg, 1932.



\bibitem{Kutateladze} 
S. S. Kutateladze, Aleksandrov of Ancient Hellas, Siberian Electronic Mathe-
matical Reports, 9 (2012), p. A6-A11.

\bibitem{KNR} S. S. Kutateladze, S. P. Novikov, and Yu. G. Reshetnyak, 
Aleksandr Danilovich Aleksandrov
(on the 100th anniversary of his birth),
Russian Math. Surveys, vol. 67 No.5  (2012)  p. 959-965



\bibitem{PT} J. H. Lambert, Theorie der Parallellinien (1766), French translation in:  La théorie des lignes parallèles de Johann
Heinrich Lambert, ed. A. Papadopoulos and G. Théret, Collection Sciences dans l'Histoire, Librairie Scientifique
et Technique Albert Blanchard, Paris, 2014.



\bibitem{Leibniz} G. W. Leibniz, Mathematische Schriften, zweite Abteilung I, Berlin 1849.


\bibitem{Papadopoulos-Minkowski} A. Papadopoulos, Busemann’s Work on the Characterization of Minkowski Geometries, In: Herbert Busemann, Selected Works, Vol. II, p. 145-159, 2018.


\bibitem{Notices} A. Papadopoulos,  Herbert Busemann, Notices of the American Mathematical Society, 2018, Vol. 65, No. 3, March 2018, p. 341-343.





 \bibitem{Papadopoulos-Hilbert} A. Papadopoulos, Hilbert's fourth problem. In: Handbook of Hilbert geometry (A. Papadopoulos and M. Troyanov, ed.), European Mathematical Society Publishing House,  IRMA Lectures in Mathematics and Theoretical Physics, Vol. 22, p. 391-431, 2014.





 \bibitem{Papa-Chronogeometry} A. Papadopoulos, Chronogeometry,  In: Herbert Busemann, Selected works, II, ed. A. Papadopoulos,
Springer, Cham, 2018
p. 133-141.


\bibitem{P-Menelaus} A. Papadopoulos, Menelaus' \emph{Spherics} in Greek and Arabic mathematics and beyond. This Handbook.

\bibitem{Pogorelov-H} A. V. Pogorelov, Hilbert's fourth problem. Translated by R. A. Silverman, edited by I. Kra in cooperation with E. Zaustinskyi, V. H. Winston and Sons, Washington DC, 1979. Original published by Nauka.


\bibitem{RP2} 
R. Rashed and  A. Papadopoulos, Menelaus'  Spherics:  Early Translation and al-M\=ah\=an\=\i /al-Haraw\=\i 's Version, De Gruyter, Series: Scientia Graeco-Arabica,  21,  2017, 890 pages.





 \bibitem{R1} R. Rashed, Angles et grandeur : d’Euclide à  Kam\={a}l al-D\={\i}n al-F\={a}ris\={\i}, Scientia Graeco-Arabica , Band 17. De Gruyter, Boston/Berlin, 2015
 \bibitem{R2} R. Rashed, Apollonius : Les Coniques, tome 1.1 : Livre I, commentaire historique et
mathématique, édition et traduction du texte arabe, de Gruyter, 2008, 666 p.
 \bibitem{R3} R. Rashed, Apollonius : Les Coniques, tome 2.2 : Livre IV, commentaire historique
et mathématique, édition et traduction du texte arabe, de Gruyter, 2009, 319 p. 
 \bibitem{R4} Apollonius : Les Coniques, tome 3 : Livre V, commentaire historique et mathématique,
édition et traduction du texte arabe, de Gruyter, 2008, 550 p.
 \bibitem{R5} R. Rashed, Apollonius : Les Coniques, tome 4 : Livres VI et VII, commentaire
historique et mathématique, édition et traduction du texte arabe, Scientia Graeco-Arabica, vol. 1.4, de Gruyter, 2009, 572 p. 
 \bibitem{R6} R. Rashed, Apollonius : Les Coniques, tome 2.1 : Livres II et III, commentaire historique
et mathématique, édition et traduction du texte arabe, de Gruyter, 2010, 682 p.
 
 



\bibitem{Reidemeister} K. Reidemester, \"Uber die singulären Randpunkte eines konvexen Körpers. Math.
Ann. 83 (1921), p. 116-118.

 \bibitem{Predire} R. Thom, Pr\'edire n'est pas expliquer,  Entretiens avec Emile No\"el, Champs, Paris,  Eshel, coll. La Question, 1991.
 
 


 \bibitem{Thom-Logos} R. Thom, Expos\'e introductif, In: Logos et Th\'eorie des Catastrophes, J. Petitot (dir.), Actes du colloque international de Cerisy de 1982, Centre culturel international de Cerisy-la-Salle,  Genève, Pati\~{n}o, 1989. Publié aussi dans  
 Gazette des math\'ematiciens,  édition sp\'eciale, 2004,  p. 7-22.

\bibitem{Zalgaller} V.  A. Zalgaller, Memoirs of A. D. Aleksandrov and his Leningrad geometry seminar, Siberian Electronic Mathematical Reports, 9 (2012), p. A26-A80.

 

  \end{thebibliography}
\end{document}